\documentclass{article}
\usepackage[14pt]{extsizes}
\usepackage[utf8]{inputenc}
\usepackage[english,russian]{babel}
\usepackage{math}
\usepackage{xcolor}
\usepackage{amssymb}
\usepackage{amsmath}
\usepackage{graphics}
\usepackage{tikz}
\usepackage{pgfplots}
\usepackage{comment}

\title{Внешние биллиарды вне правильного восьмиугольника: периодичность почти всех орбит и существование апериодической орбиты}
\author{Филипп Рухович}

\begin{document}

\maketitle

Раздел: теория динамических систем

517.938 УДК

\begin{abstract}
Доказано существование апериодической орбиты для внешнего биллиарда вне правильного восьмиугольника, а также, что почти все траектории такого внешнего биллиарда являются периодическими; явно выписаны все возможные периоды.
\end{abstract}

Пусть $\gamma$ - выпуклая фигура, а $p$ — точка вне ее. Проведем правую относительно $p$ касательную к $\gamma$; определим $Tp \equiv T(p)$ как точку, симметричную $p$ относительно точки касания.

\begin{Def} \label{base}
Отображение $T$ называется внешним биллиардом; фигура $\gamma$ называется столом внешнего биллиарда.
\end{Def}
\begin{Def}
Точку $p$ вне фигуры $\gamma$ назовем периодической, если существует такое натуральное $n$, что $T^np = p$; минимальное такое $n$ назовем периодом точки $p$ и обозначим как $per(p)$. 
\end{Def}
\begin{Def}
Точку $p$ вне фигуры $\gamma$ назовем апериодической, если она --- не периодическая, а ее траектория бесконечна в две стороны. 
\end{Def}
\begin{Def}
Точку $p$ вне фигуры $\gamma$ назовем граничной, если $T^np$ не определено для некоторого $n \in {Z}$.
\end{Def}

В данной статье будем полагать, что $\gamma$ --- выпуклый многоугольник.

Внешние биллиарды были введены Бернардом Нойманном в 1950-х годах и стали популярны в 1970-х благодаря Ю.Мозеру \cite{Moser78}. Внешние биллиарды исследовались рядом авторов\ (см. например, \cite{Tab93}, \cite{Schwartz09}, \cite{DF09}, а также монографию \cite{Tab05}). Так, Р.Шварц \cite{Schwartz09} показал, что траектория начальной точки может быть неограниченной, тем самым разрешив вопрос Мозера\ -\ Нойманна, поставленный в \cite{Moser78}.

В центре нашего внимания находятся следующие открытые в общем случае \underline{проблемы периодичности}:

\begin{enumerate}
\item Существует ли апериодическая точка для внешнего биллиарда вне правильного $n$-угольника?
\item Какова мера периодических орбит внешнего биллиарда вне правильного $n$-угольника?
\end{enumerate}

С.Л.Табачников в \cite{Tab93} решил проблемы периодичности для случая правильных 3-х, 4-х и 6-угольника: для них апериодической точки нет, а также для 5-угольников - для них апериодические точки существуют, но их мера равна нулю. В монографии \cite{Tab05}, опубликованной Американским математическим обществом в 2005 г., С.Л.Табачников приводит результаты компьютерного моделирования для восьмиугольника, но пишет, что <<строгого анализа до сих пор нет>>, и что для других случаев результатов нет. В дальнейшем правильный пятиугольник и связанная с ним символическая динамика подробно исследовались в работе N.Bedaride и J.Cassaigne \cite{BC11} (см. также их монографию \cite{BC11a}).

В монографиях \cite{Schwartz10}, \cite{Schwartz14} Р.Шварц исследовал внешний биллиард вне правильного восьмиугольника и множество связанных с ним вопросов; однако решения проблем периодичности им получено не было.

Основным результатом данной работы является следующие теоремы.
\begin{Th} \label{MainTreorem}
Для внешнего биллиарда вне правильного восьмиугольника существует апериодическая точка.
\end{Th}

\begin{Th} \label{MainTreorem2}
В случае внешнего биллиарда вне правильного восьмиугольника, периодические точки образуют вне стола множество полной меры.
\end{Th}

\begin{Th} \label{MainTreorem3}
Множество всевозможных периодов точек для внешнего биллиарда вне правильного восьмиугольника есть объединение следующих множеств:

\begin{enumerate}
\item $\{12*9^n-4*(-3)^n, 12*9^n+4*(-3)^n, 4*9^n | n \in \mathrm{Z}, n \geq 0\}$;
\item $\{8, 8k, 8k*9^n| n, k \in \mathrm{Z}, n \geq 0, k \geq 2\}$;
\item $\{24k*9^n-4*(-3)^n+12*9^n| n, k \in \mathrm{Z}, n \geq 0, k \geq 2\}$;
\item $\{24k*9^n+4*(-3)^n-4*9^n| n, k \in \mathrm{Z}, n \geq 0, k \geq 2\}$.
\end{enumerate}

\end{Th}
        
Доказательство теорем базируется на наличии самоподобных структур внутри
инвариантной относительно внешнего биллиарда вне правильного восьмиугольника фигуры. Шварц в \cite{Schwartz10} высказал гипотезу о существовании самоподобных структур в случаях 5, 10, 8, 12; однако с тех пор доказательства не появилось.

\begin{figure}[h!]
\begin{center}
\includegraphics[width=80mm]{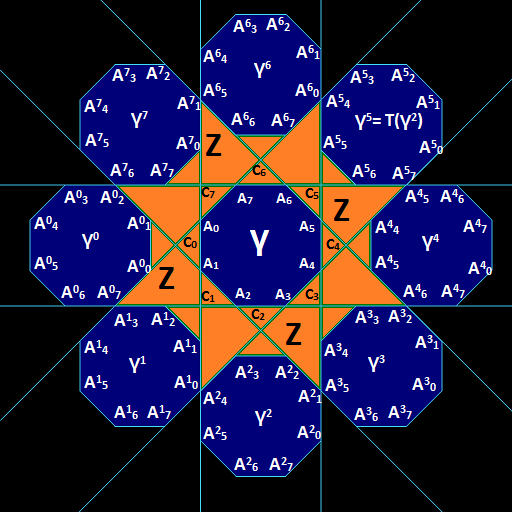}
\caption{Внешний биллиард вне правильного восьмиугольника: фигуры $\gamma$, $\gamma^i$, $Z$}.
\label{ris:picEnterZRaw}
\end{center}
\end{figure}

Опишем самоподобные структуры. Для этого рассмотрим (рис. \ref{ris:picEnterZRaw}) правильный восьмиугольник $\gamma := A_0A_1A_2A_3A_4A_5A_6A_7$, в котором вершины занумерованы против часовой стрелки, и преобразование внешнего биллиарда $T$. Пусть $l_i$ есть прямая, проходящая через вершины $A_i$ и $A_{(i+1)\ mod\ 8}$, $i = 0, 1, \ldots, 7$. Определим также точку $C_i, i\in[0,7]$ как точку пересечения прямых $l_{(i - 1)\ mod\ 8}$ и $l_{(i+1)\ mod\ 8}$. Пусть $\gamma^i = A^i_0A^i_1A^i_2A^i_3A^i_4A^i_5A^i_6A^i_7$ - это открытый многоугольник, симметричный столу $\gamma$ относительно точки $C_i$.

Многоугольники $\gamma^i$ образуют <<ожерелье>> вокруг $\gamma$ и соприкасаются с соседями по вершинкам; так, $A^1_7 = A^2_3$, $A^2_0 = A^3_4$ и т.д.. Таким образом, $\gamma^i$ ограничивают связную область - многоугольник $Z := A^0_0A^0_7A^0_6A^1_1A^1_0A^1_7 \ldots A^7_7A^7_6A^7_5$.

Сформулируем несколько очевидно следующих из рис. \ref{ris:picEnterZRaw} утверждений.
\begin{Lm}
$\forall i \in \{0,1,\ldots,7\}: T(\gamma^i) = \gamma^{(i+3)\ mod\ 8}$.
\end{Lm}

\begin{Lm}
$T(Z) \subset Z \supset T^{-1}(Z)$.
\end{Lm}

Будем исследовать точки внутри многоугольника $Z$. Заметим, что $Z$ и преобразование $T$ инвариантны относительно поворота на $k\pi/4, k \in \mathbb{Z}$. Отождествим точки $Z$ относительно такого поворота; тем самым область изучения можно ограничить до четырехугольника $A_1A^2_4A^2_3A^2_2$, на котором внешний биллиард $T$ индуцирует преобразование $T'$.

\begin{figure}[h!]
\begin{center}
\includegraphics[width=100mm]{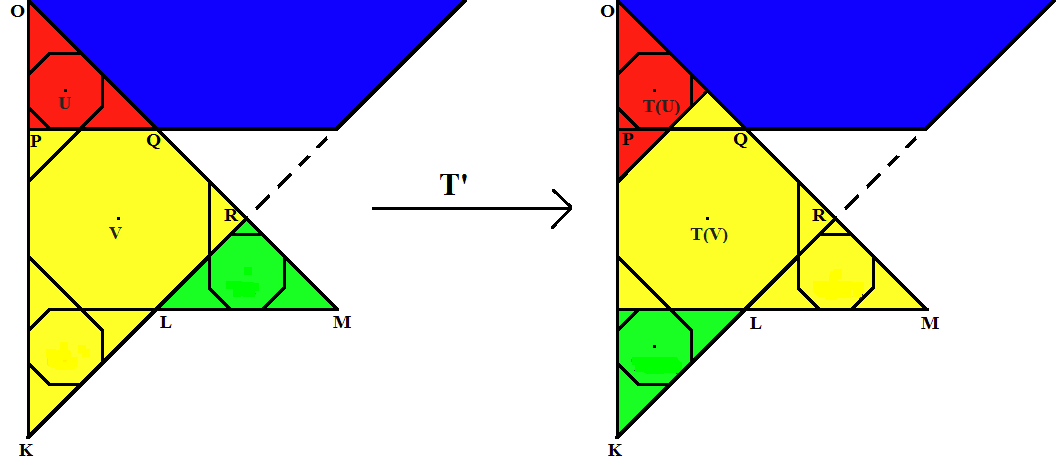}
\caption{Индуцированное преобразование $T'$ на четырехугольнике $OKLM$}.
\label{ris:picT}
\end{center}
\end{figure}

Для удобства дальнейшего анализа переобозначим точки внутри $A_1A^2_4A^2_3A^2_2$ (рис. \ref{ris:picT}). Обозначим $A_1$ как $O$, $A^2_4$ как $K$, $A^2_3$ как $L$, $A^2_2$ как $M$, $C_3$ за $R$, $A_2$ за $Q$. Пусть $P$ - точка пересечения отрезка $OK$ и луча $A_3A_2$, а $S$ - пересечение отрезка $OK$ и луча $RM$. 

Проведем лучи $A_3A_2$ и $A_4A_3$; они разобьют $OKLM$ на три многоугольника $OPQ$, $PQRK$, $LMR$; для каждого из этих многоугольников, $T$ есть центральная симметрия относительно вершин $A_2$, $A_3$ и $A_4$ соответственно. Будем называть эти многоугольники $W_1$, $W_2$ и $W_3$ соответственно.

\begin{Lm} \label{inducedT}
Индуцированное преобразование $T'$ есть:

\begin{enumerate}
\item для открытого треугольника $int(W_1) = int(OPQ)$ --- поворот на $\frac{3\pi}{4}$ вокруг точки $U$, являющейся точкой пересечения биссектрисы угла $KOM$ и серединного перпендикуляра к $A_1A_2$;

\item для открытого четырехугольника $int(W_2) = int(PQRK)$ --- поворот на $\frac{\pi}{2}$ вокруг точки $V$, являющейся точкой пересечения биссектрисы угла $KOM$ и биссектрисы угла $LRQ$;

\item для открытого треугольника $int(W_3) = int(LRM)$ --- поворот на $\frac{\pi}{4}$ вокруг точки $W$, являющейся точкой пересечения серединного перпендикуляра $LM$ и биссектрисы угла $KOM$.
\end{enumerate}

На общих же границах описанных многоугольников, т.е. $PQ$ и $KR$, $T'$ не определено.
\end{Lm}
 
\begin{Lm}
Точка $x \in OKLM$ является граничной (периодической, апериодической) точкой относительно преобразования $T'$, если и только если $x$ является граничной (периодической, апериодической) точкой относительно преобразования внешнего биллиарда $T$.
\end{Lm}

Изучим преобразование $T'$ на четырехугольнике $OKLM$.

\begin{Def}
Пусть $U_1$ есть правильный восьмиугольник со сторонами, параллельными сторонам $\gamma$, <<вписанный>> в треугольник $OPQ$.
\end{Def}

\begin{Def}
Обозначим за $\Gamma$ аффинное преобразование, т.ч. $\Gamma (O) = O$, $\Gamma(A^3) = U_1$.  Также введем точки $K' := \Gamma(K)$, $L' := \Gamma (L)$, $M' := \Gamma (M)$.
\end{Def}

Отметим, что $\Gamma(OKLM) = OK'L'M'$.

Введем преобразование $T''$, играющее ключевую роль для дальнейших исследований.

\begin{Def} \label{firstReturnMap}
Пусть $x \in int(OK'L'M')$. Пусть $n_x$ есть минимальное такое натуральное число $n$, что $T'^n(x) \in int(OK'L'M')$. Тогда $T''(x) := T'^{n_x}(x)$. Если же такого $n_x$ не существует, то $T''(x)$ не определено.
\end{Def}

\begin{Lm} \label{selfSimilarity}
Пусть $x \in int(OKLM)$. Тогда:
\begin{enumerate}
\item $T'(x)$ определено, если и только если $T''(\Gamma (x))$ определено;
\item Если $T'(x)$ определено, то $\Gamma (T'(x)) = T''(\Gamma(x))$. 
\end{enumerate}
\end{Lm}

\begin{figure}[h!]
\center{\includegraphics[width=80mm]{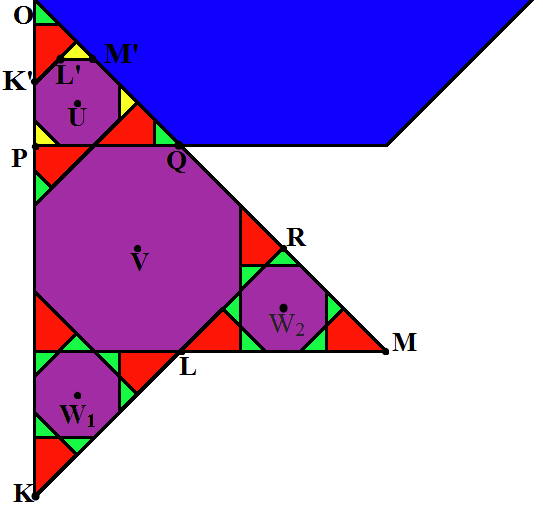}}
\caption{Траектория преобразования $T''$}.
\label{ris:picTSelf}
\end{figure}

Преобразования типа $T''$ известны в литературе как {\it преобразования первого возвращения (first return map)}. Так, $T''$ есть преобразование первого возвращения $OK'L'M'$ относительно $T'$.
Траектории первого возвращения можно увидеть на рис. \ref{ris:picTSelf}. Отметим, что множество точек, $T'$-траектории которых не проходят через $OK'L'M'$, могут быть легко описаны как четыре правильных восьмиугольника; все эти восьмиугольники суть периодические фигуры преобразования $T'$.

Поймем расположение периодических компонент в OKLM.
\begin{Def}
 Пусть $h_{-1}$ есть длина стороны восьмиугольника $\gamma^2$ (равного восьмиугольнику $\gamma$), а $h_0$ есть длина стороны восьмиугольника $U_2$. Пусть $\lambda := \frac{h_0}{h_{-1}}$. Определим также последовательность чисел $h_n := h_0 * \lambda ^ n$, $n \in \mathbb{N}$.
\end{Def}

С помощью леммы \ref{selfSimilarity} можно доказать следующую лемму.
\begin{Lm} \label{countOfFixedSize}
$\forall n \in \mathrm{Z}_+$: внутри четырехугольника $OKLM$ существует ровно $3^k$ периодических фигуры, являющихся правильными восьмиугольниками со стороной $h_n$.
\end{Lm}

\begin{figure}[h!]
\begin{center}
\includegraphics[width=80mm]{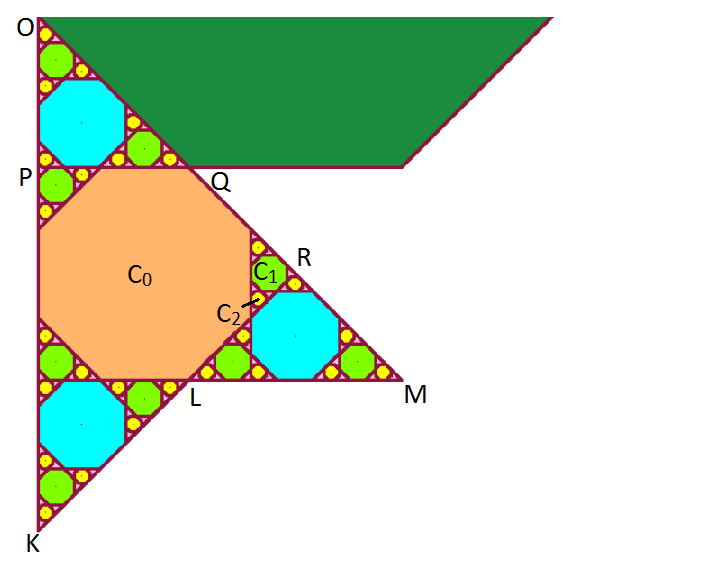}
\caption{Внешний биллиард вне правильного восьмиугольника: периодические фигуры}
\label{ris:picPeriodicCompsAndSpiral}
\end{center}
\end{figure}

Местонахождение всех периодических фигур в OKLM можно описать следующим образом. Внутри OKLM существует периодическая фигура со стороной $h_0$. Положение этой фигуры однозначно: это восьмиугольник, вписанный в четырехугольник PQRK. Удалим этот восьмиугольник из OKLM; тогда OKLM распался на три подобных OKLM с коэффициентом $\lambda$ четырехугольника, в которых по лемме \ref{countOfFixedSize} лежат три правильных восьмиугольника со стороной $h_1$. Их положения также восстанавливаются однозначно. Удалим и эти три открытых многоугольника. Оставшееся от OKLM множество точек распадается на 9 многоугольников, подобных $OKLM$ с коэффициентом $\lambda^2$; в них по лемме \ref{countOfFixedSize} находятся 9 правильных восьмиугольников со сторонами $h_2 = h_0*\lambda^2$; восстанавливаем однозначно положения этих восьмиугольников, удаляем и получаем 27 четырехугольников, подобных OKLM. Выполним эту операцию бесконечное число раз. Тогда на $k$-ом шаге, $k \leq 0$, мы получаем $3^k$ четырехугольников, подобных $OKLM$ с коэффициентом $\lambda^k$, в которых по лемме \ref{countOfFixedSize} находятся $3^k$ периодических правильных восьмиугольников со сторонами длины $h_k = h_0 * \lambda^k$, параллельными $\gamma$. Положения этих многоугольников восстанавливаются однозначно. Таким образом, периодические компоненты выглядят так, как показано на рис. \ref{ris:picPeriodicCompsAndSpiral}.

Из данного рассуждения и конфигурации периодических компонент очевидно следует

\begin{Lm} \label{MainTheoremForZ}
Периодические точки внутри фигуры $Z$ образуют множество полной меры.
\end{Lm}

Рассмотрим три периодические компоненты $C_0$, $C_1$, $C_2$, со сторонами длин $h_0$, $h_2$, $h_3$ соответственно, изображенные на рис. \ref{ris:picPeriodicCompsAndSpiral}. Они ограничивают четырехугольник $G_0$, подобный OKLM. Пусть $g$ - аффинное преобразование, переводящее OKLM в $G_0$. Пусть $C_{i+3} := g(C_i)$ и $G_{i+1} := g(G_{i}), i \geq 0$. Тогда $C_i$ есть последовательность периодических компонент, периоды которых стремятся к бесконечности, а каждый из четырехугольников $G_i$ ограничен тремя периодическими компонентами $C_{3i}$, $C_{3i+1}$, $C_{3i+2}$. Пусть точка $c$ - предел последовательности $C_n$.

\begin{Lm}
$c$ - апериодическая точка.
\end{Lm}

Следующий факт позволяет доказать теоремы 2 и 3 не только для фигуры Z, но для всей плоскости. Пусть $T_4$ есть преобразование первого возвращения $T'$ для угла $\angle A^2_0A^2_1A^3_6$. Пусть преобразование $H$ есть параллельный перенос, т.ч. $H(A_1) = A^2_1$.  

\begin{Lm} \label{selfSimilarity4}
Пусть $x \in int(\angle KOM)$. Тогда:
\begin{enumerate}
\item $T'(x)$ определено, если и только если $T_4(H(x))$ определено;
\item Если $T'(x)$ определено, то $H(T'(x)) = T_4(H(x))$. 
\end{enumerate}
\end{Lm}

Пусть $Y$ есть объединение четырехугольника $OKLM$ и восьмиугольника $\gamma^2$.

\begin{Lm} \label{allBasePeriodicComponents4}
Все периодические компоненты внутри $\angle KOM$ могут быть получены из периодических компонент, лежащих в фигуре $Y$, с помощью последовательности преобразований $H$ и $T'$. Более формально, пусть $C$ есть периодическая компонента. Тогда $C = T'^m(H^n(C_0))$, где $C_0$ - лежащая в $Y$ периодическая компонента, $n, m$ - некоторые целые неотрицательные числа, причем $m \leq 3$.
\end{Lm}

О результатах данной статьи были сделаны доклады на Combinatorics on Words, Calculability and Automata research school (CIRM, Marseille, France, 3 February 2017), а также на Tiling Dynamical System research school (CIRM, Marseille, France, 21 November 2017) .

Работа поддержана грантом РНФ № 17-11-01337.

\bibliographystyle{utf8gost705u}
\bibliography{biblio}

\end{document}